%% file: Hopkins_Levitzki.tex
\documentclass[12pt]{amsart}
\usepackage{amsmath,amssymb,ifthen}

\input{defs}

\title[Hyperradical - Hopkins--Levitzki Theorem]{The Hyperradical and The Hopkins--Levitzki Theorem for Modular Lattices}
\author{Fernando Guzm\'an
}

\newcounter{amasrefsORbibtex}
\newcommand{\amsrefsORbibtex}[2]{\ifthenelse{\isodd{\value{amasrefsORbibtex}}}{#1}{#2}}

\amsrefsORbibtex{\usepackage[initials,author-year]{amsrefs}}%
{\newcommand{\cites}[1]{\cite{#1}}}

\begin{document}

\begin{abstract}

  Many arguments in the Theory of Rings and Modules are, on close
  inspection, purely Lattice theoretic arguments.  C\v{a}lag\v{a}reanu
  \nocite{calagareanu} has a long repertoire of such results in his
  book.  The Hopkins-Levitzki Theorem is interesting from this point
  of view, because a special case of it lends to an obvious lattice
  theory approach, but the rest is a little more subtle.  Albu and
  Smith have obtained some sufficient conditions for
  the question of when Artinian implies Noetherian. Here we present a
  new approach, using the concept of Hyperradical; we obtain necessary
  and sufficient conditions.
\end{abstract}

\maketitle

\section{Introduction}


A ring is left (resp. right) Artinian, if its lattice of left
(resp. right) ideals satisfies the descending chain condition.  It is
left (resp. right) Noetherian, if its lattice of left (resp. right)
ideals satisfies the ascending chain condition.  A classical theorem
connecting these concepts is the Hopkins-Levitzki Theorem
\cites{hopkins,levitzki} (HLT for short).  It states that every left
(resp. right) Artinian ring is left (resp. right) Noetherian.  The
statement of this theorem is lattice theoretic and it is only natural
to ask if there is lattice theoretic proof of it, i.e. if it can be
extended to lattices, and under what assumptions.  The goal of this
paper is to answer that question.

The standard proof of the of the Hopkins-Levitzki Theorem found in
algebra textbooks like \cites{jacobson2, grillet} has two components.
The first one considers the special case when the Jacobson radical,
defined as the intersection of all maximal ideals, is trivial.  This
component of the HLT readily extends to lattices; see
Proposition~\ref{prop:hopkins,levitzki}. In the second component, the ring
operations play an essential role, via the nilpotency of the Jacobson
radical.  So, our question reduces to how this part of the proof can be
extended to lattices. Albu and Smith
\cites{albu,albu-smith-1,albu-smith-2} have obtained some results
related to our question.  In \cites{albu} the lattice is assumed to be
modular and upper continuous, a condition weaker than algebraic; but a rather
technical additional hypothesis, condition $(\lambda)$, is needed.  In
\cites{albu-smith-1}, the hypothesis of upper continuity is further weakened
to condition $\mathcal{E}$, which ensures that the lattice $L$ has a
good supply of essential elements, and condition $\mathcal{BL}$ which
places a bound on the composition length of some subintervals.  All
three conditions, $(\lambda)$, $\mathcal{E}$, and $\mathcal{BL}$ are
local and existential.

Since the lattice of left (right) ideals of a ring is modular
and algebraic, it makes sense to being in the lattice case with these
two assumptions. It can be seen that modularity alone will not do, by
considering $(\N,|)$, the lattice of natural
numbers under divisibility.  This lattice is Artinian but not
Noetherian.  Even though this lattice is distributive, hence modular,
it is not algebraic.  Thus the question arises if every Artinian
modular algebraic lattice is Noetherian.  The answer is no as we will
show in Example~\ref{exam:Z p infinity}.

The concept leading to the solution of our problem is the
hyperradical, which is a global construction.  As we will show, being
hyperradical free is not only a sufficient condition, but also
necessary for a modular Artinian lattice to be Noetherian.  There is
no need to assume the lattice to be algebraic, or even upper
continuous. 

It should be noted that for the lattice of left ideals of a ring, the
hyperradical free condition is not vacuous, as will be shown in
Example~\ref{exam:germs}. 

\section{The Radical, Modularity and Chain Conditions}

In Ring Theory, there are a number of different {\sl radical}
constructions.  Some of them, but not all, can be expressed as the
intersection of maximal objects in some lattice.  One of the best
known examples is the Jacobson radical of a ring $R$, $J(R)$, which is
equal to the intersection of all maximal two-sided ideals.  It is also
equal to the intersection of all maximal left (right) ideals.  The
extension of this definition to a complete lattice is immediate.

\begin{definition}\label{def:radical}
  In a complete lattice $L$, the \defem{radical} of $L$, $\rad(L)$ is
  the meet of all coatoms of $L$.  If $A$ is a (universal) algebra, we
  denote by $\rad(A)$ the radical of the lattice $\Sub(A)$ of
  subalgebras of $A$.
\end{definition}

So, the Jacobson radical $J(R)$ is the radical in the lattice of
two-sided ideals, as well as the radical in the lattice of left
ideals, and in the lattice of right ideals.  In other words, it is the
radical of $R$, when viewed as a left $R$-module.

It is easy to check that in any lattice $L$, for any $x,y,z\in L$
\[ x \leq z \imp x\join(y\meet z) \leq (x\join y)\meet z \]

\begin{definition}\label{def:modular lattice}
A lattice $L$ is \defem{modular} if it satisfies 
\[ x \leq z \imp x\join(y\meet z) = (x\join y)\meet z \]
for any $x,y,z\in L$.
\end{definition}

If $M$ is a left-module over a ring $R$, the lattice $\sub(M)$ of
submodules of $M$ is a modular lattice. In particular, the
lattice of left-ideals of $R$ is modular.

\begin{definition}\label{def:chain condition}
  A lattice $L$ is said to be \defem{Noetherian} or to satisfy the
  \defem{ascending chain condition}, ACC for short, if it contains no
  infinite  ascending chain $x_0<x_1<x_2<\cdots$\\ 
  Dually, $L$ is said to be \defem{Artinian} or to satisfy the
  \defem{descending chain condition}, DCC for short, if it contains no
  infinite descending chain $x_0>x_1>x_2>\cdots$
\end{definition}

A \defem{left-Noetherian Ring} is a ring $R$ such that the lattice of
left ideals has the ACC.  A \defem{left-Artinian Ring} is a ring $R$
such that the lattice of left ideals has the DCC.  We can rephrase the
HLT as follows: Let $R$ be a ring and $L$ the lattice of left ideals
of $R$. If $L$ if Artinian, then $L$ is Noetherian.

Not every complete Artinian lattice is Noetherian, as illustrated by
the lattice $(\N,|)$.  Even though this lattice is distributive, hence
modular, it is not algebraic.  Since the lattice $L$ of left ideals of
a ring $R$ is modular and algebraic, we ask the following question:

\begin{question}\label{q:Is Amal Noetherian}
Is every Artinian modular algebraic lattice, Noetherian?
\end{question}

The next example shows that the answer is {\bf no}, and we need to
modify the hypotheses.

\begin{example}\label{exam:Z p infinity}
Let $L$ be the
lattice of subgroups of $Z_{p^\infty}$. $L$ is a chain isomorphic to
$\langle\N\cup\{\infty\},\leq\rangle$, so it is modular and Artinian,
but it is not Noetherian.  
\end{example}

A look at the radical of $L$ in Example~\ref{exam:Z p infinity} gives us a clue of what goes wrong in
this example.  $L$ has no coatoms, hence $\rad(L)=1$.  This never
happens in the lattice of ideals of a ring (with 1), where maximal
proper ideals are always guaranteed to exist, and therefore 
$\rad(L) < 1$.

\ignore{
As it turns out, the {\sl algebraic} condition on $L$ plays no role,
but a restriction on the radical series, as defined below, is
crucial. 
}

The following lemma illustrates one of the key features of modular
lattices, which we can loosely describe as follows: ``Behavior in the
lattice can be moved inside intervals without much loss''.

\begin{lemma}\label{lemma:chains in modular lattice}
Let $L$ be a modular lattice, $a\in L$, and $(x_i)_{i\in I}$ a chain
in $L$.  Let $y_i=a\meet x_i$ and $z_i=a\join x_i$.  If the chain 
$(x_i)_{i\in I}$ is an infinite ascending (resp. descending) chain,
then so is at least one of $(y_i)_{i\in I}$ and $(z_i)_{i\in I}$.
\end{lemma}

\proof Let's consider the ascending case.  The descending case is
dual.\\
Suppose $(x_i)_{i\in I}$ is infinite ascending, but both
$(a\meet x_i)_{i\in I}$ and 
$(a\join x_i)_{i\in I}$ become stationary at 
$u=a\meet x_k=a\meet x_{k+1}=\dots$ and 
$v=a\join x_k=a\join x_{k+1}=\dots$.\\
Using modularity we get
\[\begin{array}{rcl}
x_k &=& x_k\join u \\
    &=& x_k\join(a\meet x_{k+1}) \\
    &=& (x_k\join a)\meet x_{k+1} \\
    &=& v\meet x_{k+1} \\
    &=& x_{k+1}
\end{array} \]
contradicting the assumption about $(x_i)_{i\in I}$.
\endproof

\begin{corollary}\label{cor:noetherian segments}
  Let $L$ be a modular lattice, $a\in L$.  $L$ is Noetherian
  (resp. Artinian) iff $[0,a]$ and $[a,1]$ are Noetherian
  (resp. Artinian).
\end{corollary}

\begin{corollary}\label{cor:chains under coatoms}
  Let $L$ be a modular lattice, and $m\in L$ a coatom. 
  If $(x_i)_{i\in I}$ is an infinite ascending chain in $L$, then so
  is $y_i=m\meet x_i$. 
\end{corollary}

The usual proof of the HLT, first considers
the case when the radical is $0$, and then the general case.  The
argument in the first case is lattice theoretic, as we illustrate
next.  Note that there is no assumption about the lattice being
algebraic. 

\begin{proposition}\label{prop:hopkins,levitzki}
  Let $L$ be a complete modular lattice. If $L$ is Artinian and
  radical free, i.e. $\rad(L)=0$, then $L$ is Noetherian.
\end{proposition}

\proof Being Artinian, $\rad(L)$ must be expressible as the meet of
finitely many maximal elements $m_1,\dots,m_k$.  If we had an infinite
ascending chain $(x_i)_{i\in I}$, repeated application of
Corollary~\ref{cor:chains under coatoms} yields an infinite ascending
chain $(x_i\meet m_1\meet \cdots \meet m_k)_{i\in I}$.  But $x_i\meet
m_1\meet \cdots \meet m_k \leq m_1\meet \cdots \meet m_k =\rad(L)=0$.
\endproof

This proof shows that if something is going to go wrong about $L$
being Noetherian, it
will show up below the radical.  So, we look at the interval
$[0,\rad(L)]$, and the radical of this lattice.  This gives rise to
the radical series.

\begin{definition}\label{def:loewy radical series}
Let $L$ be a complete lattice.  We define the \defem{Loewy radical
  series} of $L$ as follows: 
\begin{itemize}
  \item $r_0(L)=1$,
  \item for any ordinal $\sigma$, $r_{\sigma +
    1}(L)=\rad([0,r_\sigma(L)])$,
  \item for a limit ordinal $\sigma$, $\ds r_\sigma(L)=\bigmeet_{\alpha <
    \sigma}r_{\alpha}(L)$. 
\end{itemize}

The smallest ordinal $\sigma$ such that $r_{\sigma+1}(L)=r_\sigma(L)$
is called the \emph{Loewy radical length} of $L$, and $r_\sigma(L)$ is
called the \defem{hyper-radical} of $L$.  It is denoted by
$r_\infty(L)$.  We say that $L$ is \defem{hyper-radical free} if 
$r_\infty(L)=0$.
\end{definition}

Being hyper-radical free is precisely the extra condition needed to
extend the HLT to complete modular lattices.

\begin{theorem}\label{thm:hyper-radical free}
Let $L$ be a complete modular lattice.  If $L$ is Artinian and hyper-radical
free, then $L$ is Noetherian.
\end{theorem}

\proof Being Artinian, $L$ must have finite Loewy radical length.
Therefore, $r_n(L)=0$ for some $n\in\N$.  For $i=1,\dots,n$, the
interval $[r_i(L),r_{i-1}(L)]$, is modular; it is Artinian by
Corollary~\ref{cor:noetherian segments}; it is radical free by
construction. By Proposition~\ref{prop:hopkins,levitzki} it is Noetherian.
By Corollary~\ref{cor:noetherian segments}, $L$ is Noetherian.
\endproof

The Loewy radical series is the dual construction of the \emph{Loewy
  (socle) series}, see \cites{calagareanu}.
The dual of radical-free is semiatomic, i.e. when
the socle is equal to $1$.  The dual of the hyper-radical we call the
\emph{hyper-socle}.  A lattice is \emph{hyper-semiatomic} if the
hyper-socle is equal to $1$.  By duality we get the following theorem:

\begin{dual}\label{thm:full hyper-socle}
Let $L$ be a complete modular lattice.  If $L$ is Noetherian and
hyper-semiatomic, then $L$ is Artinian. 
\end{dual}

After extending the HLT to hyper-radical free modular lattices, a
number of questions arise.  Is the ``hyper-radical free'' hypothesis
necessary?  Is it vacuous?  We answer these questions with an example
and a proposition.  But first a lemma.

\begin{lemma}\label{lemma:radical of module}
Let $R$ be a ring, $N$ a left $R$-module, and $J$ the (Jacobson)
radical of $R$.  Then $J\cdot N\leq\rad(N)\leq N$
\end{lemma}

\proof If $M$ is a maximal submodule of $N$ then $N/M$ is simple and
$J$ is contained in $\ann(N/M)$.  In other words, $J\cdot N\leq M$.
It follows that $J\cdot N\leq\rad(N)$.
\endproof

\begin{example}\label{exam:germs}
  The hyper-radical free hypothesis in Theorem~\ref{thm:hyper-radical
    free} is not vacuous, not even for the ideal lattice of a ring.
  The ideal lattice \mbox{$L=\Idl(R)$} of the ring $R$ of germs of
  $C^{\infty}(\R)$ functions at $x=0$, is not hyper-radical free, and
  its Loewy radical length is $\omega$. To see this, note that $R$ is
  a local ring with maximal ideal 
  \[ M=\{f\in R|f(0)=0\}=x\cdot R, \] 
  so this is the first radical of $L$, i.e. $J=M$.  $M$ has a single
  maximal submodule 
  \[ M_2=\{f\in M|f'(0)=0\}=x^2\cdot R, \] 
  so this is the second radical of $L$, and inductively,
  \[ r_n(L)=\{f\in R|f^{(i)}(0)=0\text{ for } i=0,\dots,n-1\}=x^n\cdot
  R. \]  
  Therefore, 
  \[ r_{\omega}(L)=\ds\bigintersection_n \left(x^n\cdot
  R\right)= \{f\in R|f^{(i)}(0)=0\text{ for all } i\}, \] 
  the ideal of germs of flat functions.  Now, by
  Lemma~\ref{lemma:radical of module} we have 
  \mbox{$J\cdot r_{\omega}(L)\leq r_{\omega+1}(L)\leq 
  r_{\omega}(L)$}.  But $J\cdot r_{\omega}(L)= x\cdot R\cdot
  r_{\omega}(L)= x\cdot r_{\omega}(L)= r_{\omega}(L)$, so
  $r_{\omega+1}(L)=r_{\omega}(L)$.   It is a well-known fact that
  there are non-zero flat functions, like $f(x)=\exp(-1/x^2)$. So,
  $r_{\infty}(L)=r_{\omega}(L)\neq 0$.
\end{example}

The hyper-radical free hypothesis in Theorem~\ref{thm:hyper-radical
  free} is necessary.

\begin{proposition}\label{prop:noetherian is hyper-radical free}
  If $L$ is a complete Noetherian lattice, then it is hyper-radical
  free.
\end{proposition}

\begin{proof}
  If we had $r_\infty(L)>0$, then the interval $[0,r_\infty(L)]$ would
  have no coatoms.  Therefore, it must have an infinite ascending
  chain.
\end{proof}

This proposition, combined with the Hopkins-Levitzki Theorem, tells us
that the lattice of left ideals of a left-Artinian ring is
hyper-radical free.  The proof given of Theorem~\ref{thm:hyper-radical
  free}, does not replace the standard proof of the HLT for rings,
unless one finds a direct argument to show that the lattice of left
ideals of a left-Artinian ring is hyper-radical free.

Example~\ref{exam:germs} was suggested by Mazur and
Karagueuzian \cites{mazur-karagueuzian}.  It would be nice to have a
characterization of the class of rings which are hyper-radical free.

\amsrefsORbibtex{

}
{
\bibliography{radical}
\bibliographystyle{plain}
}

\noindent
Department of Mathematical Sciences\\
State University of New York at Binghamton\\
Binghamton, NY 13902-6000 \\
fer@math.binghamton.edu\\

\end{document}

%% file: defs.tex
\expandafter\ifx\csname theoremstyle\endcsname\relax\else%
\theoremstyle{plain}\fi

\newtheorem{theorem}{Theorem}
\newtheorem{xdual}{Theorem}
\newtheorem{corollary}[theorem]{Corollary}
\newtheorem{lemma}[theorem]{Lemma}
\newtheorem{proposition}[theorem]{Proposition}
\newtheorem{question}[theorem]{Question}

\expandafter\ifx\csname theoremstyle\endcsname\relax\else%
\theoremstyle{definition}\fi

\newtheorem{definition}[theorem]{Definition}
\newtheorem{example}[theorem]{Example}

\newenvironment{dual}{%
\begin{xdual}}
{\end{xdual}}

\newcounter{alfa}

\newcommand{\ignore}[1]{\relax}

\newcommand{\ds}{\displaystyle}
\newcommand{\join}{\vee}

\newcommand{\meet}{\wedge}
\newcommand{\bigmeet}{\bigwedge}

\newcommand{\bigintersection}{\bigcap}

\newcommand{\imp}{\Rightarrow}
\newcommand{\sub}{{\rm Sub}}
\newcommand{\Sub}{\textrm{Sub}}
\newcommand{\Idl}{\textrm{Idl}}

\newcommand{\rad}{{\rm r}}
\newcommand{\N}{{\mathbb N}}

\newcommand{\R}{{\mathbb R}}

\newcommand{\ann}{\text{Ann}}

\newcommand{\defem}[1]{%
  \ifmmode{\ #1\ }\else{{\mdseries\slshape\sffamily #1}}\fi\index{#1}}